\documentclass[12pt]{article}

\usepackage{amssymb,amsmath,latexsym,mathrsfs}

\usepackage{url}
\usepackage{color}
\usepackage{authblk}
\usepackage{comment}
\usepackage{graphicx}
\usepackage{subfigure}
\usepackage{caption}
\usepackage{algorithm}
\usepackage[noend]{algpseudocode}
\usepackage{cite}
\usepackage[margin=3cm]{geometry}

\newcommand{\btheta}{\boldsymbol{\theta}}

\newcommand\bX{\mathbf{X}}

\usepackage{amsthm}

\begin{document}

\title{\LARGE \bf Parameters estimation in a 3-parameters $p$-star model}

\author{Pietro Lenarda$^1$, Giorgio Gnecco$^1$, Massimo Riccaboni$^1$ \\
\small{$^1$ IMT, School for Advanced Studies, Lucca, Italy}}
\vspace{5mm}

\maketitle

\thispagestyle{plain}
\pagestyle{plain}


\begin{abstract}
An important issue in social network analysis refers to the development of algorithms
for estimating optimal parameters of a social network model, using data available from the network itself.
This entails solving an optimization problem. In the paper, we propose a new method for parameters estimation in a specific social network model, namely, the so-called $p$-star model with three parameters. The method is based on the mean-field approximation of the moments associated with the three subgraphs defining the model, namely: the mean numbers of edges, two-stars, and triangles. A modified gradient ascent method is applied to maximize the log-likelihood function of the $p$-star model, in which the components
of the gradient are computed using approximate values of the moments. Compared to other existing iterative methods for parameters estimation, which are computationally very expensive when the number of vertices becomes large, such as gradient ascent applied to maximum log-pseudo-likelihood estimation, the proposed approach has the advantage of a much cheaper cost per iteration,
which is practically independent of the number of vertices.
\end{abstract}

\section{Introduction}

\subsection{\bf Motivation of the work.} Since the 1970s, researchers in several fields have used social network analysis
to investigate interpersonal social relationships, communication networks, scientific paper co-authorships
and citations, patterns in protein interaction, and so on \cite{Goldenbergetal2009}. Nowadays, this issue is of particular relevance due to the presence of online networking communities such as Facebook and LinkedIn. An important aspect of social network analysis refers to the development of algorithms for estimating optimal parameters of a social network model, using data available from the network itself. This entails solving an optimization problem, such as maximum likelihood estimation. Unfortunately, applying algorithms such as (exact) gradient ascent to solve the optimization problem is often unfeasible, due to high computational cost needed to evaluate the gradient exactly. For this reason, in the paper, a modified gradient ascent method, based on the so-called mean-field approach\footnote{Loosely speaking, mean-field theory investigates the behaviour of complex stochastic models composed of a large number of mutually interacting units, by reducing them to simpler models, in which the effect on any given unit of all the other units is approximated by a single averaged effect, and fluctuations around such an average effect are neglected. Mean-field approximations have been used in optimization (especially in combinatorial optimization \cite{orland85} and variational inference \cite{Xingetal2002}), e.g., in connection with optimization algorithms such as simulated annealing \cite{bilbroetal92}.}, is proposed for maximum likelihood estimation, for a particular social network model, known as the $p$-star model, focusing on its important case with three parameters.

\subsection{\bf Preliminaries}. The $p$-star model, also called exponential random graph model, is one of the best-known and widely used statistical models for social networks \cite{Goldenbergetal2009,goodreau2007,lusher2012exponential,robins2007intro,snijders2011statistical,Snijdersetal2006,Handcock,Wainwright}. The model has been applied with success in applications related to fields such as communication, computer science, physics, psychology, and sociology \cite{ShuPal2010}. Compared with previous social network models, the main reason for its success is that the $p$-star model is able to represent interdependencies in a network. The term ``$p$-star'' was coined by Wasserman and Pattison \cite{WassermanPattison} in honor of the work \cite{HollandLeinhardt} by Holland and Leinhardt, which appear to be the first in the literature to have proposed a specific one-parameter instance of the model, called $p_1$ therein. A more general version of the model in \cite{HollandLeinhardt}, based on a Markovianity assumption, was proposed by Frank and Strauss in \cite{FrankStrauss1986}, and called Markov graph model therein\footnote{Such Markov $p$-star model can be extended to the case in which the Markovianity assumption does not hold \cite{PattisonRobins2002}. The class of $p$-star models includes both cases. Since the Markov case is more frequent in applications, in the paper we often omit the term ``Markov'' when referring to the Markov $p$-star model by Frank and Strauss.}. The sufficient statistics of this model (called $\mathbb{P}_{\btheta}$ in the following) are expressed in terms of subgraphs called $p$-stars and triangles. In the model, the parameters vector $\btheta$ contains the parameters associated with each among the subgraphs above, and provides information about the macroscopic properties of the social network.

In principle, the problem of estimating the parameters vector $\btheta$ of the $p$-star model $\mathbb{P}_{\btheta}$ is solved by applying iterative methods for the maximization of the log-likelihood function obtained from the real data. In doing this, one needs to use, in each iteration, exact values of the moments associated with the three subgraphs defining the
model, namely: the mean numbers of edges, two-stars, and triangles. However, because of the high computational cost needed to compute the so-called partition function associated with the $p$-star model, an exact evaluation of such moments (and as a consequence, the estimation problem itself) becomes intractable when the number of vertices $n$ is large.
To circumvent this drawback, instead of working directly with the log-likelihood function, the log-pseudo-likelihood function 
\cite{Handcock} is also used in the literature. This function is defined in terms of vectors of differences of statistics
$\Delta \bX^k_{ij}$, one for each sample $\bX^k$ and each edge $(i,j)$. Nevertheless, also the computation of these vectors requires a high computational cost when the number of vertices $n$ of the network is large.

The 
mean-field approach to estimate the moments of a $p$-star model was proposed originally by Park and Newman 
\cite{ParkNewman}, limiting to the case of a $2$-dimensional parameters vector, with parameters associated, respectively, with the edges and the triangles. After that work, Chatterjee and Diaconis \cite{Diaconis} solved the model asymptotically for the non-negative high-dimensional case, i.e., for a finite-dimensional parameters vector in which all the parameters (apart from the first one) are non-negative, and the number of vertices tends to $+\infty$.

{\bf Contributions of the work.} In this paper, following the work \cite{ParkNewman} by Park and Newman, we derive a mean-field approximation to compute approximately the moments of a $3$-parameters $p$-star model. 
In this model, the features are the numbers of edges, 2-stars, and triangles. Moreover, we apply the mean-field approximation to estimate the parameters vector of such a $p$-star model.
Roughly speaking, denoting by $\mathbb{E}[E(\bX)]$, $\mathbb{E}[S_2(\bX)]$ and $\mathbb{E}[T(\bX)]$, respectively, the 
mean numbers of edges (i.e., 1-stars), 2-stars, and triangles in a network generated under a $p$-star model $\mathbb{P}_{\btheta}$ in which the parameters $p$, $q$ and  $r$ are, respectively, the ``average probabilities'' that an edge, a $2$-star and a triangle is present in the network, 
the mean-field approach consists in computing approximately
\begin{eqnarray}
 \mathbb{E}[E(\bX)] &=& p \binom{n}{2}, \nonumber \\
 \mathbb{E}[S_2(\bX)] &=& n\binom{n-1}{2} q, \nonumber \\
 \mathbb{E}[T(\bX)] &=& \binom{n}{3}  r, \nonumber
\end{eqnarray}
using suitable approximations of $p$, $q$, and $r$, where $\binom{n}{2}$, $n \binom{n-1}{2}$, and $\binom{n}{3}$, are the maximum possible numbers of edges, $2$-stars and triangles in a network with $n$ vertices. 
To find the mean-field approximations for $p$, $q$ and $r$, we derive appropriate equations linking $p$, $q$, and $r$, which we then solve numerically. The justification for using the mean-field approximation is that when the parameters vector $\btheta$ is in the so-called high-temperature phase
\cite{Mixingtime}, the $p$-star model $\mathbb{P}_{\btheta}$ behaves like an independent Erd\H{o}s-R\'{e}nyi model $G(n, p^*)$ \cite{Bollobas}, for a given $p^*$.

The proposed method is tested against maximum log-pseudo-likelihood estimation, confirming the computational advantage
of the mean-field approximation, which requires only the solution of a nonlinear system of equations having practically constant
computational cost per iteration (i.e., independent of the number of vertices $n$, which can be interpreted as a measure of the size of the problem).

{\bf Organization of the work.} The paper is structured as follows. In Section \ref{pstar_model}, the $p$-star model $\mathbb{P}_{\btheta}$ for social networks is described, focusing on the case of a model with $3$ parameters. Sections \ref{sec:ML} and \ref{sec:MPL} deal, respectively, with maximum log-likelihood and maximum log-pseudo-likelihood estimation of the parameters vector in a $3$-parameters $p$-star model. In Section \ref{sec:MF}, the proposed mean-field approximation of some moments in the $3$-parameters $p$-star model is introduced, and applied to maximum log-likelihood estimation of the parameters vector. In Section \ref{sec:MH}, the Metropolis-Hastings sampler is described, for the generation of the samples used to define the log-likelihood and log-pseudo-likelihood functions. Section \ref{sec:examples} compares the proposed parameters estimation method with maximum log-pseudo-likelihood estimation, demonstrating its computational advantages considering both a simulated and a real social network. Finally, Section \ref{sec:conclusions} concludes the paper with a discussion of our results and future research directions. Details about the implementation of Newton's method for the approximate solution of the system of nonlinear equations derived from the mean-field approximation are reported in the Appendix.
 
\section{The $p$-star model $\mathbb{P}_{\boldsymbol{\theta}}$ for social networks}\label{pstar_model}
 In the present context, a social network 
of $n$ individuals is represented by a symmetric random matrix, in which each entry $\bX_{ij}$ is a random variable accounting for the 
type of relation between the individuals $i$ and $j$. In the simplest case $\bX_{ij} \in \{0,1 \}$ assumes only binary values (i.e., $\bX$ is an adjacency matrix), and $X_{ii}=0$ for each $i$.
Let $\bX:=(X_{ij})$ be an $n \times n$ symmetric zero-one matrix, representing a graph in which the presence or not of an edge is a binary random
variable $X_{ij}$ (self-loops are excluded). Under the so-called Markovianity assumption, i.e., assuming that $X_{ij}$ (with $i \neq j$) and $X_{kl}$ (with $k \neq l$) are conditionally independent
given the rest of $\bX$ if and only 
if $\{i,j\} \cap \{k,l\}= \emptyset$, the probability distribution (Gibbs measure) over the set of all symmetric zero-one matrices with zero diagonal can be parametrized as follows \cite{FrankStrauss1986}:
\begin{equation} \label{eq:gen}
 \mathbb{P}_{\boldsymbol{\theta}}(\bX=X):= \exp \left( \sum_{i,j,k} \theta_{ijk} \mathbf{1}_{T_{ijk}}(X)+
 \sum^n_{k=1} \sum_{i_0, \dots , i_k} \theta_{i_0 \dots i_k} \mathbf{1}_{S_{i_0 \dots , i_k}}(X)- A(\btheta) \right),
\end{equation}
where $\mathbf{1}_{T_{ijk}}(X)$ and $\mathbf{1}_{S_{i_0 \dots , i_k}}(X)$ are indicator functions of triangles and $k$-stars in the network, and $A(\btheta)$ is a normalizing factor, which makes the sum of the probabilities of all the possible realizations of the matrix be equal to 1. A triangle is defined as a set of three edges $T_{ijk}:=\{ (i,j), (j,k), (k,i) \}$ with $i\neq j \neq k$, and a $k$-star is a set of edges $S_{i_0, \dots , i_k}:=\{ (i_0, i_1), (i_0, i_2), \dots , (i_0, i_k) \}$ in which $i_0 \neq i_1 \neq \ldots, i_k$, and one of the vertices is always $i_0$. Notice that the $1$-stars are the edges.
The proof of Equation \eqref{eq:gen} above is a consequence of the well-known Hammersley-Clifford theorem \cite{Besag,HammCliff,Lauritzen}. The model (\ref{eq:gen}) is referred to in the paper as the (Markov) $p$-star model (or $\mathbb{P}_{\boldsymbol{\theta}}$). 

Under an additional homogeneity assumption, and truncating the expansion above to the order two, we can rewrite Equation  \eqref{eq:gen}  
as the following 3-parameters $p$-star model:
\begin{equation}\label{eq:3params}
 \mathbb{P}_{\boldsymbol{\theta}}(\bX=X):=\exp \left( \theta_1 E(X)+ \theta_2 S_2(X)+ \theta_3 T(X)- A(\btheta)  \right),
\end{equation}
where $E(X)$, $S_2(X)$, $T(X)$ are, respectively, the numbers of edges, 2-stars, and triangles in the network $X$, and
$\boldsymbol{\theta}:=(\theta_1, \theta_2, \theta_3)$ is the parameters vector. Qualitative properties of the network are associated with the 3 parameters above as follows:
\begin{enumerate}
\item the edges parameter $\theta_1$ can be regarded as a measure of the density of the network;
\item the $2$-stars parameter $\theta_2$ is a measure of the tendency of the network to clustering;
\item the triangles parameter $\theta_3$ is a measure of transitivity.
\end{enumerate} 
The $3$-parameters $p$-star model is defined as $\mathbb{P}_{\boldsymbol{\theta}}(\bX=X)=\exp \left( -H(X) - A(\btheta) \right)$, where $H(X)$ is the Hamiltonian, defined as follows:
\begin{align}\label{eq:hamiltonian}
 H(\bX)&:=\theta   \sum_{i <j} X_{ij}+ \sigma  \sum_{k \neq i,j} \sum_{i <j} X_{ik} X_{ki}- \alpha \sum_{i <j <k} X_{ij} X_{jk} X_{ki} \nonumber \\
       &=\theta E(\bX)+ \sigma S_2(\bX)-\alpha T(\bX),
\end{align}
where $A(\btheta):=\log \left( \sum_{X} e^{-H(X)} \right)$ is the normalizing factor (called log-partition function, since $Z:=\exp(A(\btheta))$ is the partition function), and
the parameters vector is also defined as
$\boldsymbol{\theta}:=(\theta_1, \theta_2, \theta_3):=(-\theta, -\sigma, \alpha)$. 

\section{Maximum log-likelihood estimation}\label{sec:ML}

We address now the problem of estimating $\btheta$ on the basis of observed data.
Consider a $3$-parameters $p$-star model $\mathbb{P}_{\btheta}$ as in Equation \eqref{eq:3params}. 
Let $\{ \bX^1, \dots , \bX^N  \}$ be a collection of independent and identically distributed (i.i.d.)
networks sampled from $\mathbb{P}_{\btheta}$. The aim is to maximize the log-likelihood function, defined as follows:
\[ L(\btheta; \bX^1, \dots , \bX^N ):= \frac{1}{N} \sum^N_{i=1} \log{\mathbb{P}_{\btheta}}( \bX^i) .\]
Let $\hat{\boldsymbol{\mu}}:=(\hat{\mu}_1, \hat{\mu}_2, \hat{\mu}_3)$ be the vector of empirical moments, i.e.,
the empirical mean numbers of edges, 2-stars, and triangles obtained from the data:
\[ \hat{\mu}_1:=\frac{1}{N} \sum^N_{i=1} \mathbf{1}_{E}(\bX^i) , \quad
   \hat{\mu}_2:=\frac{1}{N} \sum^N_{i=1} \mathbf{1}_{S_2}(\bX^i) , \quad 
   \hat{\mu}_3:=\frac{1}{N} \sum^N_{i=1} \mathbf{1}_{T}(\bX^i), \]
where $\mathbf{1}_{E}(\bX), \mathbf{1}_{S_2}(\bX), \mathbf{1}_{T}(\bX)$ are the numbers of edges, 2-stars, and triangles in the network $X$.
The log-likelihood can be also written as:
\begin{equation}\label{eq:male}
L(\btheta; \bX^1, \dots , \bX^N  )= \btheta \cdot \hat{ \boldsymbol{\mu} } - A(\btheta).
\end{equation}
The Maximum Log-Likelihood Estimate (MLLE) is the vector of parameters $ \btheta^*_{MLLE}$ maximizing the objective function \eqref{eq:male}. 
It can be shown \cite{Wainwright} that \eqref{eq:male}, as a function of $\btheta $, is a concave function, and that
its maximum exists. The gradient of $L(\btheta; \bX^1, \dots , \bX^N )$ is given by:
\begin{equation}\nonumber
\nabla_{\btheta} L(\btheta)= \hat{ \boldsymbol{\mu} } - \boldsymbol{\mu}\,,
\end{equation}
where $\boldsymbol{\mu}:=(\mu_1, \mu_2, \mu_3)$ 
is the vector whose components are the expected numbers of edges, 2-stars, and triangles under the $3$-parameters $p$-star model $\mathbb{P}_{\btheta}$, which are defined as follows:
\begin{equation}\label{eq:expectations}
\mu_1:=\mathbb{E}[E(\bX)],\quad \mu_2:=\mathbb{E}[S_2(\bX)], \quad \mu_3:= \mathbb{E}[T(\bX)].
\end{equation}
In the following, to find the maximum of the log-likelihood function, the gradient ascent method \cite{Bertsekas} is applied. The iterative algorithm reads as follows:
 \begin{algorithm}[H]
\caption{ \textbf{ Maximum Log-Likelihood Estimation (MLLE) via gradient ascent} } \label{alg1}
\begin{algorithmic}
\State \text{Fix a stepsize $\gamma >0$, a number of iterations $N_{\rm it}$, and an initialization $\btheta^{(0)}$ for $\btheta$.}
\State \text{Evaluate the empirical moments vector $\hat{\boldsymbol{\mu}}$.}
\For {$k=0,1,\dots,N_{\rm it}-1$} 
\State Evaluate the moments vector: 
$$
\boldsymbol{\mu}^{(k)}=\boldsymbol{\mu}^{(k)}(\btheta^{(k)}).
$$
\State \text{Update the parameters vector as}:
\begin{align*} 
\btheta^{(k+1)}&=\btheta^{(k)}+ \gamma\left(\hat{\boldsymbol{\mu} }- \boldsymbol{\mu} ^{(k)} \right).
\end{align*}
\EndFor
\end{algorithmic}
\end{algorithm}
Of course, also variations of the algorithm above can be used in principle, including, e.g., the case of a variable stepsize, the insertion of a different termination criterion, and the application of coordinate maximization and coordinate gradient ascent \cite{BecTet2013,LiuWright2015,Nesterov2012,Wright2015}. Even if the algorithm described above looks very appealing to solve the maximum log-likelihood estimation problem, there is a serious issue related to the computation of the moments vector $\boldsymbol{\mu}^{(k)}$ for each parameters vector
$\btheta^{(k)}$ generated by the algorithm. Indeed, the exact evaluation of the moments becomes rapidly nearly intractable as the number of vertices $n$ grows, due, to the computational difficulties 
related, respectively, to the need of considering all the possible realizations of the matrix $\bX$ when computing the expected value in (\ref{eq:expectations}), and of evaluating the log-partition function $A(\btheta^{(k)})$ (which is also defined in terms of all such possible realizations). In the literature, there exist iterative methods for estimating the exact moments, such as as the
elimination algorithm and the junction tree algorithm \cite{Wainwright},
but they work well only for small $n$. In this paper, we propose a different approach based on the so-called mean-field approximation of the moments vector $\boldsymbol{\mu}$. Before introducing such an approach, in the next section we describe another strategy for estimating the parameters vector, which is based on the maximization
of the log-pseudo-likelihood function \cite{Handcock}. Then, we compare the two approaches.

\section{Maximum log-pseudo-likelihood estimation}\label{sec:MPL}

In this section, the problem of maximizing the log-pseudo-likelihood function is addressed. This method is widely used in parameters estimation for general exponential random graph models, because iterative methods for maximum log-pseudo-likelihood estimation are typically computationally more tractable when compared to those for maximum log-likelihood estimation. Under appropriate assumptions, it is known
(see, e.g., \cite{Koller}) that, when the number of samples $N$ tends to $+\infty$, the Maximum Log-Pseudo-Likelihood Estimate (MLPLE) of the parameters vector converges to its maximum log-likelihood estimate.

For a $p$-star model $\mathbb{P}_{\btheta}$ and an edge $(i,j)$, we set $p_{ij}:=\mathbb{P}_{\btheta}(X_{ij}=1 \vert \bX_{-ij})$,
where $\bX_{-ij}$ is the collection of all the remaining edges. Let $\btheta:=(\theta_1, \theta_2, \theta_3)^T$ be the
parameters vector, and $E(X), S_2(X), T(X)$ be the corresponding numbers of $1$-stars, $2$-stars, and triangles for a network $X$. Then, for each edge $(i,j)$, the vector $\Delta \bX_{ij}$ of difference statistics is defined as follows: 
$$ \Delta \bX_{ij}:=( E(X^+_{ij})-E(X^-_{ij}),
                   S_2(X^+_{ij})-S_2(X^-_{ij}),
                     T(X^+_{ij})-T(X^-_{ij}) )^T,
                          $$
where $X^{+}_{ij}$ is the network associated with the matrix with $X_{ij}=1$ and all the remaining entries equal to the corresponding entries of $X$, and $X^{-}_{ij}$ is the network associated with the matrix with $X_{ij}=0$ and all the remaining entries equal to the corresponding entries of $X$. The log-pseudo-likelihood function associated with the model $\mathbb{P}_{\btheta}$ and the samples $\bX^1, \dots , \bX^N$ is defined as follows:
\begin{align*}
PL(\btheta; \bX^1, \dots , \bX^N)&:=
\dfrac{1}{N} \sum^{N}_{k=1} \sum_{i,j} \log \mathbb{P}_{\btheta}(X^k_{ij}=1 \vert \bX^k_{-ij}) \\
&=\dfrac{1}{N} \sum^{N}_{k=1} \sum_{i,j} \left[ X^{k}_{ij} \log(p^k_{ij})+ (1- y^{k}_{ij}) \log(1-p^k_{ij})  \right]\\
&=\dfrac{1}{N} \sum^{N}_{k=1} \sum_{i,j} \left[ X^{k}_{ij} \btheta \cdot \Delta \bX^k_{ij} - \log(1+ \exp( \btheta \cdot  \Delta \bX^k_{ij}) ) \right].
\end{align*}
One can notice that the expression above is equivalent to the log-likelihood for a logistic regression model, in which each element of the 
adjacency matrices $X^1_{ij}, \dots, X^N_{ij}$ is treated as an independent observation, with the corresponding row of the 
design matrix given by $\Delta \bX^1_{ij}, \dots, \Delta \bX^N_{ij}$.

Since the function ${PL}(\btheta; \bX^1, \dots , \bX^N)$ is concave in the parameters vector
$\btheta$, we apply the gradient ascent
method to approximate the optimal $\btheta^*_{MLPLE}$, where the gradient of the objective function has the following expression:
\begin{equation*}
\nabla_{\btheta}  PL(\btheta^{(k)})= \dfrac{1}{N} \sum^{N}_{k=1} \sum_{i,j} \left( X^{k}_{ij}  - \dfrac{ \exp( \btheta^{(k) } \cdot  \Delta \bX^k_{ij}  ) }{1+  \exp( \btheta^{(k) } \cdot   \Delta \bX^k_{ij}  ) } 
\right) \Delta \bX^k_{ij}.
\end{equation*}
The gradient ascent method, applied to the maximization of the log-pseudo-likelihood function, reads as follows:
 \begin{algorithm}[H]
\caption{ \textbf{ Maximum Log-Pseudo-Likelihood Estimation (MLPLE) via gradient ascent} } \label{alg2}
\begin{algorithmic}
\State \text{Fix a stepsize $\gamma >0$, a number of iterations $N_{\rm it}$, and an initialization $\btheta^{(0)}$ for $\btheta$.}
\State \text{Evaluate the empirical moments vector $\hat{\boldsymbol{\mu}}$.}
\For {$k=0,1,\dots,N_{\rm it}-1$} 
\State \text{Update the parameters vector as}:
\begin{align*} 
\btheta^{(k+1)}&=\btheta^{(k)}+ \gamma \nabla_{\btheta}  PL(\btheta^{(k)}).
\end{align*}
\EndFor
\end{algorithmic}
\end{algorithm}

\section{Mean-field approximation of the moments}\label{sec:MF}

In this section, we derive explicit formulas for the approximate computation of the moments
of the features associated with a 3-parameters $p$-star model, using the mean-field approach. The starting point of our analysis is the paper \cite{ParkNewman}, in which the problem was addressed for a 2-parameters $p$-star model.
Consider all the terms in the Hamiltonian \eqref{eq:hamiltonian} involving the edge $X_{ij}$, for $i \neq j$:
\begin{align*}
 \partial H(X_{ij}=x):= x \left[\theta+ \sigma \sum_{l \neq i,j} X_{il}+ \sigma \sum_{l \neq i,j} X_{jl}- \alpha \sum_{l \neq i,j} X_{jl} X_{li} \right].
\end{align*}
Then, the mean number of edges $p:= \mathbb{E}[X_{ij}]$ is: 
\begin{align*}
 \mathbb{E} [X_{ij}]
 &=\sum_{X} \mathbb{P}(X) \dfrac{\mathbb{P}(X_{ij}=1)}{\mathbb{P}(X_{ij}=1)+\mathbb{P}(X_{ij}=0)} \\
 &=\dfrac{1}{Z} \sum_{X} \mathbb{P}(X) \dfrac{e^{-\partial H(X_{ij}=1)}}{e^{-\partial H(X_{ij}=1)}+e^{-\partial H(X_{ij}=0)}}\\
 &=\mathbb{E} \left[ \dfrac{1}{e^{\left( \theta +\sigma \sum_{l \neq i,j} X_{il}+\sigma \sum_{l \neq i,j} X_{jl} -\alpha \sum_{l \neq i,j}X_{jl} X_{li} \right) }+1} \right].
\end{align*}
Let $q$ be the mean number of 2-stars under the $3$-parameters $p$-star model. Then, approximating in the previous expression all the terms of the form
$X_{jl}X_{li}$ ($j \neq l \neq i$) with $q$, we obtain the following approximate expression for $p$:
\begin{equation}\label{eq:p}
 p \approx \dfrac{1}{e^{\theta -\alpha(n-2)q+2\sigma (n-2)p}+1}.
\end{equation}
Now, let
\begin{align*}
 \partial H(X_{ij}=x,X_{jk}=y):=&\theta x+ \theta y
                          -\alpha x \sum_{l \neq i,j,k} X_{jl} X_{li} 
                          -\alpha y \sum_{l \neq i,j,k} X_{lj} X_{lk} 
                          -\alpha xy X_{ki} \\
                          &+\sigma x \sum_{l \neq i,j} X_{il}
                           +\sigma y \sum_{l \neq j,k} X_{kl}
                           +\sigma (x+ y) \sum_{l \neq i, j,k} X_{jl}
                           +\sigma xy
 \end{align*}
 be the collection of all the terms in the Hamiltonian involving $X_{ij}$ and $X_{jk}$ ($i \neq j \neq k$). Then, the mean number of 2-stars $q= \mathbb{E} [X_{ij}X_{jk}]$ is expressed as follows:

{
\begin{align*}
&\sum_{X} \mathbb{P}(X) \dfrac{\mathbb{P}(X_{ij}=1, X_{jk}=1)}{\mathbb{P}(X_{ij}=1, X_{jk}=1)+\mathbb{P}(X_{ij}=1, X_{jk}=0)+\mathbb{P}(X_{ij}=0, X_{jk}=1)+\mathbb{P}(X_{ij}=0, X_{jk}=0)} \\
&=\dfrac{1}{Z} \sum_{X} \mathbb{P}(X) \dfrac{e^{-\partial H(X_{ij}=1, X_{jk}=1)}}{e^{-\partial H(X_{ij}=1, X_{jk}=1)}
                                            +e^{-\partial H(X_{ij}=0, X_{jk}=1)}
                                            +e^{-\partial H(X_{ij}=1, X_{jk}=0)}
                                            +e^{-\partial H(X_{ij}=0, X_{jk}=0)}}\\
&=\mathbb{E} \left[ \dfrac{e^{(\alpha X_{ki}-\sigma)}}{D}
                    \right],
\end{align*}
}

where
\begin{align*}
D:= &\left( e^{(\theta-\alpha \sum_{l \neq i,j,k} X_{jl} X_{li}+ \sigma \sum_{l \neq i,j} X_{il}+\sigma \sum_{l \neq i, j,k} X_{jl})}+1 \right) \\
   & \times \left( e^{(\theta-\alpha \sum_{l \neq i,j,k} X_{lj} X_{lk}+ \sigma \sum_{l \neq j,k} X_{kl}+\sigma \sum_{l \neq i, j,k} X_{kl})} +1 \right) 
    +(e^{(\alpha X_{ki}-\sigma)}-1).
                          \end{align*}

Now, because $e^{\alpha X_{ki}}$=$1 + (e^{\alpha} -1)X_{ij}$, passing to the mean-field approximation, we obtain the following approximate expression for $q$:
\begin{equation}\label{eq:q}
 q \approx \dfrac{e^{- \sigma}(1+(e^{\alpha} -1)p)}{\left( e^{\theta-\alpha (n-3)q +\sigma (2n-5) p}+1 \right)^2+ e^{- \sigma}(1+(e^{\alpha} - 1)p)-1 }.
\end{equation}
Finally, for $i \neq j \neq k$, consider
\begin{align*}
&\partial H(X_{ij}=x,X_{jk}=y, X_{ki}=z)\nonumber \\
:=& -\alpha x \sum_{l \neq i,j,k} X_{il} X_{lj} -\alpha y \sum_{l \neq i,j,k} X_{jl} X_{lk} 
     -\alpha z \sum_{l \neq i,j,k} X_{kl} X_{li} 
     -\alpha xyz \\
 & + \sigma (x+z) \sum_{l \neq i,j,k} X_{il}
+\sigma (x+y) \sum_{l \neq i,j,k} X_{jl} 
+\sigma (y+z) \sum_{l \neq i,j,k} X_{kl} \\
& +\sigma xy
+\sigma xz
+\sigma yz+  \theta x+ \theta y + \theta z.
 \end{align*}


Then, the mean number of triangles $r=\mathbb{E}[X_{ij} X_{jk} X_{ki}]$ is expressed as follows:
\begin{equation}\nonumber
 \mathbb{E}[X_{ij} X_{jk} X_{ki}]=\mathbb{E}\left[ \dfrac{e^{\alpha}}{\Delta+(e^{\alpha}-1)} \right],
\end{equation}
where:
\begin{align*}
 \Delta:=& \left( e^{(\theta -\alpha \sum_{l \neq i,j,k } X_{jl} X_{lj}+\sigma \sum_{l \neq i,j,k} X_{il}+ \sigma \sum_{l \neq i,j,k} X_{kl} + \sigma)} +1\right) \\
        & \times \left( e^{(\theta -\alpha \sum_{l \neq i,j,k } X_{jl} X_{lk}+\sigma \sum_{l \neq i,j,k} X_{jl}+ \sigma \sum_{l \neq i,j,k} X_{kl} + \sigma)} +1\right) \\
        & \times \left( e^{(\theta -\alpha \sum_{l \neq i,j,k } X_{jl} X_{li}+\sigma \sum_{l \neq i,j,k} X_{il}+ \sigma \sum_{l \neq i,j,k} X_{jl} + \sigma)} +1\right).
\end{align*}
Using again the mean-field approximation, this leads to the following approximate expression for $r$:
\begin{equation}\label{eq:r}
 r \approx \dfrac{e^{\alpha}}{\left( e^{\theta - \alpha q(n-3)+2 \sigma (n-3)p+\sigma} +1\right)^3+ (e^{\alpha}-1)}.
\end{equation}
When replacing the approximation sign with the equality sign, Equations \eqref{eq:p} and \eqref{eq:q} form a nonlinear system of equations in the unknowns $p,q$, which we
solve numerically using Newton's method \cite{bonnansetal2006} (details are reported in the Appendix). Once the values of $p$ and $q$
are determined, they are used in \eqref{eq:r} to find $r$.
Then, the corresponding average numbers of $1$-stars, $2$-stars, and triangles are expressed as follows:
\begin{align*}
\mu_{1}  \approx \binom{n}{2}  p , \quad
 \mu_{2}  \approx n \binom{n-1}{2}  q,  \quad
  \mu_3  \approx \binom{n}{3} r. 
\end{align*}
The expressions above are approximations, since $p$, $q$ and $r$ have been computed by solving the nonlinear system derived from the mean-field approximation.
For $\sigma=0$, we obtain the same result as in \cite{ParkNewman}. It is worth noticing that the mean-field approximation works well if the number 
of vertices $n$ is large, and the components $\theta_1=-\theta, \theta_2=-\sigma, \theta_3=\alpha=$ of the parameters vector $\btheta$ are in a suitable range (called high-temperature phase), as explained in details in the next paragraphs.

Concluding, the modified gradient ascent method for the maximization of the log-likelihood function using the mean-field approximation for the computation of the moments reads as follows:

 \begin{algorithm}[H]
\caption{ \textbf{Maximum Log-Likelihood Estimation (MLLE) via modified gradient ascent based on the mean-field approximation of the moments} } \label{alg3}
\begin{algorithmic}
\State \text{Fix a stepsize $\gamma >0$, a number of iterations $N_{\rm it}$, and an initialization $\btheta^{(0)}$ for $\btheta$.}
\State \text{Evaluate the empirical moments vector $\hat{\boldsymbol{\mu}}$.}
\For {$k=0,1,\dots,N_{\rm it}-1$} 
\State Evaluate the mean-field approximation of the moments vector: 
$
\boldsymbol{\mu}^{(k)}_{\rm MF}.
$
\State \text{Update the parameters vector as}:
\begin{align*} 
\btheta^{(k+1)}&=\btheta^{(k)}+ \gamma\left(\hat{\boldsymbol{\mu} }- \boldsymbol{\mu} ^{(k)}_{\rm MF} \right).
\end{align*}
\EndFor
\end{algorithmic}
\end{algorithm}

When is one allowed to use the mean-field approximation? In 
\cite{Mixingtime}, the authors make a distinction between a high- and a low-temperature phase
for a $p$-star model with parameters vector $\btheta$. 
Considering a $3$-parameters $p$-star model of the form \eqref{eq:3params}, following \cite{Mixingtime}, we define:
\begin{equation}\nonumber
 \varphi_{\btheta}(p):=\dfrac{1}{1+\exp\left( \Psi_{\btheta}(p) \right)},
\end{equation}
where 
$\Psi_{\btheta}(p):=\theta+2 \sigma(n-2)p -\alpha(n-2)p^2$.
\begin{itemize}
 \item {\bf High-temperature phase}. We say that $\btheta$ is in the high-temperature phase if $\varphi_{\btheta}(p)=p$ has a unique solution $p^*$ such that $0<\varphi'_{\btheta}(p^*)=p^* <1$.
 \item {\bf Low-temperature phase}. We say that $\btheta$ is in the low-temperature phase if $\varphi_{\btheta}(p)=p$ has at least two fixed points $p^*$ which satisfy $0<\varphi'_{\btheta}(p^*)=p^* <1$.
\end{itemize}

One can notice that the only difference between the espression $\varphi_{\btheta}(p)$ and the espression
on the left-hand side of Equation \eqref{eq:p} is the presence of the term $p^2$ instead of $q$.
This is not surprising because, if $p$ is the
average edge connectivity in a $p$-star model, then, making a rough estimate, the average probabilities of having a 2-star and a triangle ($q$ and $r$, respectively), are approximately $q \approx p^2$ and $r \approx p^3$ \cite{Diaconis}. 

In \cite{Mixingtime}, it is shown that, when $n$ is large and for $\btheta$ in the high-temperature phase, the $p$-star model
is not appreciably different from the classical Erd\H{o}s-R\'{e}nyi random graph model. Indeed, in such case, most of the Gibbs measure of the 
model is concentrated
on configurations which are essentially indistinguishable from the ones obtained through an Erd\H{o}s-R\'{e}nyi $G(n, p^* )$
random graph model, for a given $p^*$. In this model, the edges are chosen independently.
This result justifies the use of the mean-field approximation \cite{KapWie2001} to compute
the moments of the statistics defining the model, namely the expected numbers of edges, $2$-stars, and triangles, for
a $p$-star model with parameters vector $\btheta$ in the high-temperature phase and with large $n$.
Instead, in the low-temperature phase, the mean-field approximation does not work \cite{Mixingtime}.
This is also justified by the fact that, when the parameters vector is
in the low-temperature phase, the matrix in \eqref{eq:Newton} in the Appendix
is typically close to be singular, as observed in some preliminary numerical tests.

\begin{figure}[H]
\centering
\subfigure[High-temperature phase]{\includegraphics[width=.48\textwidth]{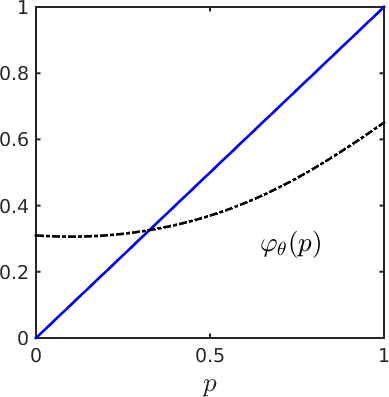}} \qquad
\subfigure[Low-temperature phase]{\includegraphics[width=.48\textwidth]{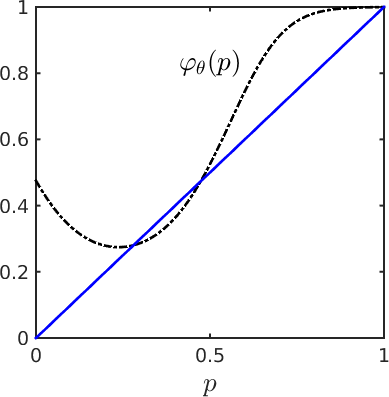}} 
\caption{$(a)$ High-temperature phase, obtained setting $\theta=-1.6, \sigma=-0.2/n, \alpha=2/n$ with $n=18$; 
         $(b)$ Low-temperature phase, obtained setting $\theta=-0.1, \sigma=-0.23, \alpha=0.97$ with $n=18$.}
\label{fig.AR1}
\end{figure}

Figure \ref{fig.AR1} refers to two parameters vectors, respectively in the high- and in the low-temperature phase.

In the 3-parameters $p$-star model considered so far, the parameters $\theta_1$, $\theta_2$, and $\theta_3$ are responsible for
the formation, respectively, of edges, 2-stars, and triangles. If the transitivity parameter $\theta_3$ is too large with respect
to the other two parameters, then the model is likely to be in the low-temperature phase, as evinced by preliminary numerical tests. Hence, a good rule of thumbs for the applicability of the method at hand (i.e., having the parameters vector in the high-temperature phase) is that the condition above is violated. This suggests to use the proposed method for estimating networks having a low number of triangles.


\section{The Metropolis-Hastings sampler}\label{sec:MH}
 
The Metropolis-Hastings algorithm is often used to obtain samples of i.i.d. graphs drawn from a $p$-star model $\mathbb{P}_{\btheta}$
\cite{Snijders2002}. In our context, this algorithm is useful to test the algorithms of the previous sections by a) fixing a parameters vector for the $p$-star model, b) generating $N$ i.i.d. samples according to that choice of the parameters vector, then c) estimating the parameters vector itself using each of the algorithms, and d) assessing the quality of each estimate (see the next Subsection \ref{subsec:example} for such a test).  
For a rigorous and exaustive analyisis of the properties of the Metropolis-Hastings algorithm when applied to exponential random graph models such as the Markov $p$-star model, we refer to \cite{Mixingtime}.
The algorithm consists in generating a Markov chain
$( \bX^t )_{t \in \mathbb{N}}$, having stationary probability distribution equal to $\mathbb{P}_{\btheta}$.
Let $\tau_{\rm mix}$ be the mixing time of the 
Markov chain. For $t \geq \tau_{\rm mix}$, the Markov chain can be assumed to be stationary, making the procedure able to produce $N$ samples $\bX^1, \dots , \bX^N$ of networks, which are essentially drawn from $\mathbb{P}_{\btheta}$, as described briefly in the following. Let $\btheta$ be a given vector of parameters for the $p$-star model $\mathbb{P}_{\btheta}$,
$t_{ \rm burn} \geq \tau_{\rm mix}$ the burn-in time of the algorithm,
$N$ the number of samples generated by the algorithm, and $m$ a given integer. The algorithm proceeds as follows: 
\begin{enumerate}
\item A random graph $\bX^0$ is generated to initialize the algorithm.
\item At each time $t \in \{0, 1, \dots , t_{\text{burn}}+ mN -1\}$:
\begin{itemize}
\item an edge $(i,j)$ in the current network $\bX^t$ is randomly chosen,
 and a new network $\bX$ is obtained from $\bX^t$ switching $(i,j)$ (from $0$ to $1$, or vice-versa).
 \item The acceptance probability 
$ \alpha(\bX^t, \bX):= \text{min}\left\{ 1, \dfrac{\mathbb{P}_{\btheta}(\bX)}{\mathbb{P}_{\btheta}(\bX^t)} \right\} $
is evaluated. 
\item A value $U$ is generated from the uniform $(0,1)$ distribution. If $U \leq \alpha ( \bX^t , \bX )$,
then set $\bX^{t + 1} = \bX$, else set $\bX^{t + 1} = \bX^t$.
\end{itemize}
\end{enumerate}
The first $t_{\rm burn}$ networks are neglected (hence the name burn-in time for $t_{\rm burn}$), then a network is selected every $m$ samples, for a total of $N$ selected networks.

For a $p$-star model with parameters vector $\btheta$ in the low-temperature phase,
the Markov chain is known to take an exponentially long time to converge to its stationary probability distribution.
In this case, sampling from the $p$-star model using the Metropolis-Hastings algorithm is highly inefficient. On the other hand,
for a $p$-star model with parameters vector $\btheta$ in the high-temperature phase, the mixing time of the Markov chain is of order $\Theta(n^2 \log n)$, and the Metropolis-Hastings algorithm is known to produce samples representative of the underlying probability distribution in a reasonably small number of steps.


\section{Numerical tests}\label{sec:examples}
We consider in what follows two examples, one with synthetic data and the other one with real data taken from the specialized 
literature. 
\subsection{Example 1}\label{subsec:example}
In order to assess the capability of the mean-field approximation for its use in the modified gradient ascent
method for maximum log-likelihood estimation,
we perform some numerical tests, by varying the number of vertices $n$, and comparing the 
computational cost of the proposed approach (Algorithm \ref{alg3}) with the one of the gradient ascent method applied to the maximization of the log-pseudo-likelihood function (Algorithm \ref{alg2}). Algorithm \ref{alg1} is not considered since it is only of theoretical interest, as explained in Section \ref{sec:ML}.

\begin{figure}[H]
\centering
\includegraphics[trim=0cm 0cm 0cm 0cm, clip=true, totalheight=0.35\textheight, angle=0]{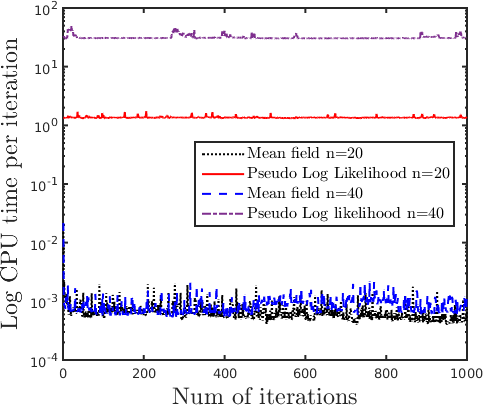}
\caption{CPU time per iteration for Algorithms \ref{alg2} and \ref{alg3} applied to the numerical tests in Subsection \ref{subsec:example}, for the two cases $n=20, 40$.}
\label{fig.CPU}
\end{figure}

The tests are performed to show that the proposed mean-field approximation of the components of
the gradient in the modified gradient ascent method for maximum log-likelihood estimation has the advantage of being computationally much faster with respect to the gradient ascent method applied to the maximization of the log-pseudo-likelihood function.
This is particularly useful when the number of vertices $n$ is large. \\
The tests are performed as follows. First, for each $n$, the 3 parameters of the vector $\btheta$ are fixed as follows:
$$ \theta_1=2 \beta_1, \quad \theta_2=\dfrac{\beta_2}{n}, \quad \theta_3=\dfrac{\beta_3}{n},
$$ with 
$\beta_1=-0.8, \beta_2=-0.2, \beta_3=2$.
Notice that the parameters for the 2-stars ($\theta_2$) and of the triangles ($\theta_3$) are 
rescaled by the factor $n$ as in \cite{Diaconis}. This is also justified by the form of the function $\Psi_{\btheta}(p)$. 
Then, a number $N=950$ of i.i.d. samples
$\{ \bX^1, \dots, \bX^N \}$, drawn from $\mathbb{P}_{\btheta}$, is generated
using the Metropolis-Hastings algorithm \cite{Snijders2002}, as described in the previous section. The number of samples $N$ is 
large enough to ensure that the log-pseudo-likelihood estimate practically coincides with the log-likelihood estimate. 
Then, the gradient ascent method is applied to the log-pseudo-likelihood function (Algorithm \ref{alg2}), and is compared to the mean-field approximation of the moments in the modified gradient ascent method applied to the log-likelihood function (Algorithm \ref{alg3}), for an increasing number of vertices $n=10, 20, \dots, 80 $, using for each comparison the same data set $\{ \bX^1, \dots, \bX^N \}$.
In Figure \ref{fig.CPU}, the CPU time per iteration is shown for the two cases $n=20, 40$. It is clear that the proposed Algorithm \ref{alg3}, based on the mean-field 
approximation of the moments, requires a much cheaper cost per iteration than Algorithm \ref{alg2}, based on the log-pseudo-likelihood function, and that its cost is practically the same for increasing values of $n$.


In contrast, when the dimension $n$ of the problem increases, the computational cost per iteration of the gradient ascent method for the maximum log-pseudo-likelihood estimation increases significantly, because the evaluation of the quantities $\Delta \bX^k_{ij}$ becomes more and more expensive. On the other hand, our proposed mean-field approximation of the moments, even if it requires more iterations for convergence to the ``true'' parameters vector (i.e., the one used to generate the samples), has the advantage of having a 
very cheap cost per iteration, because only the solution of a nonlinear system of equations is required for each iteration, which is obtained via a nested Newton's iterative procedure.

\begin{figure}[H]
\centering
\subfigure[]{\includegraphics[trim=0cm 0cm 0cm 0cm, clip=true, totalheight=0.27\textheight, angle=0]{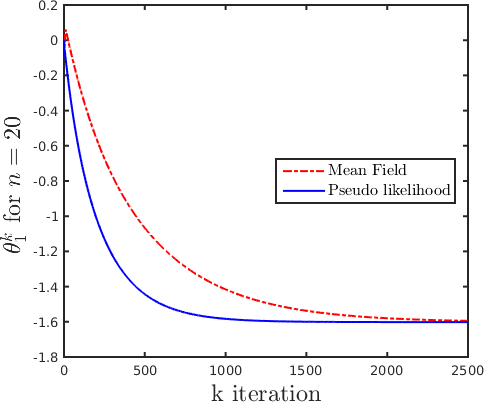}}\\
\subfigure[]{\includegraphics[trim=0cm 0cm 0cm 0cm, clip=true, totalheight=0.27\textheight, angle=0]{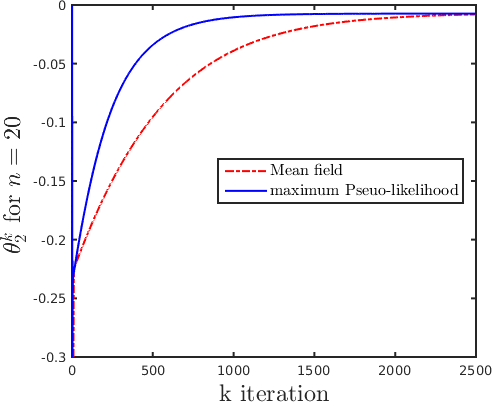}}\\
\subfigure[]{\includegraphics[trim=0cm 0cm 0cm 0cm, clip=true, totalheight=0.27\textheight, angle=0]{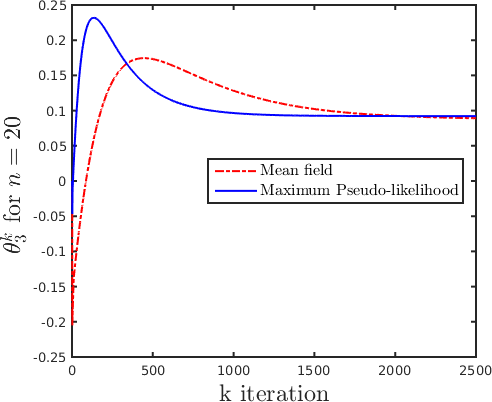}} 
\caption{Sequences of estimates of the parameters obtained by Algorithms \ref{alg2} and
\ref{alg3} for the numerical test in Subsection \ref{subsec:example} with $n=20$.}
\label{fig.Conv}
\end{figure}

Figure \ref{fig.Conv} reports, for $n=20$ and $i=1,2,3$, the sequences $\{\theta^k_i\}^{N_{\rm it}}_{k=1}$ of estimates for each parameter, obtained, respectively, when maximizing the log-likelihood via the modified gradient ascent method based on the mean-field approximation (Algorithm \ref{alg3}) and maximizing the log-pseudo-likelihood via gradient ascent (Algorithm \ref{alg2}). One can notice from the figure that the 
size of the sample $N$ is large enough to ensure that the maximum log-pseudo-likelihood estimate practically coincides with the maximum log-likelihood estimate, because the two sequences converge to the same estimates $\theta^*_i$, $i=1,2,3$. For a fair comparison, the number of iterations $N_{\rm it}$ and the stepsize $\gamma$ of the gradient ascent method is the same for the two methods in all the simulations. It is clear that, although the sequences produced by Algorithm \ref{alg2} are usually more precise than the ones produced by Algorithm \ref{alg3} in the approximation of the parameters, they require a much larger computational cost per iteration. For this reason, the computational time (and not the number of iterations) being the same, the proposed Algorithm \ref{alg3} produces better estimates than Algorithm \ref{alg2}.


Finally, Table \ref{tab:1} reports, for each simulation and $i=1,2,3$, the final estimates $\theta^*_{i, {\rm MF}}$ obtained by the proposed Algorithm \ref{alg3}, together with the corresponding components of the ``true'' parameters vector which has been used to generate the samples through the Metropolis-Hastings algorithm. The table shows, for these cases, the effectiveness of the Metropolis-Hastings sampler for the generation of the samples, since the estimates above are quite near the corresponding ``true'' parameters. It is also worth mentioning that, in all cases, the final estimated parameters are in the high-temperature phase.

\begin{table}
\caption{For the numerical tests of Subsection \ref{subsec:example}: final estimates $\theta^*_{i, {\rm MF}}$ obtained by Algorithm \ref{alg3}, and corresponding components of the ``true'' parameters vector, used to generate the samples through the Metropolis-Hastings algorithm.}
\centering
{
\begin{tabular}{ccc ccc ccc}
\hline
\noalign{\smallskip}
$n$    & $N_{\rm it} $  & $\gamma$  & $\theta_1=2\beta_1$  &$\theta^*_{1, {\rm MF}}$  & $\theta_2=\beta_2/n$ & $\theta^*_{2, {\rm MF}}$ &  $\theta_3= \beta_3/n$ &      $\theta^*_{3, \rm{MF}}$\\
\noalign{\smallskip}
\hline
\hline\noalign{\smallskip}
 \hline
\noalign{\smallskip}
10 & 1000   & 1E-2   &-1.6    &-1.5377  & -0.02     & -0.0462 & 0.2     & 0.2105 \\      
20 & 2500   & 1E-3   &-1.6    &-1.5947  & -0.01     & -0.080  & 0.1     & 0.0891 \\
30 & 50000  & 1E-4   &-1.6    &-1.5884  & -0.006667 & -0.0053 & 0.06667 & 0.0460 \\
40 & 50000  & 1E-4   &-1.6    &-1.5877  & -0.005    & -0.0060 & 0.05    & 0.0518 \\
50 & 150000 & 1E-5   &-1.6    &-1.5977  & -0.004    & -0.0036 & 0.04    & 0.0355 \\
60 & 150000 & 1E-5   &-1.6    &-1.5873  & -0.003333 & -0.0036 & 0.03333 & 0.0296 \\
70 & 150000 & 1E-5   &-1.6    &-1.5775  & -0.002857 & -0.0032 & 0.02857 & 0.0250 \\ 
80 & 900000 & 1E-6   &-1.6    &-1.6054  & -0.0025   & -0.0025 & 0.0250   &  0.0266 \\
\noalign{\smallskip}\hline
\end{tabular}
}
\label{tab:1}
\end{table}

\subsection{Example 2: The trade network of Renaissance Florentine families}
\begin{figure}[H]
\centering
{\includegraphics[trim=1cm 1cm 1cm 1cm, clip=true, totalheight=0.29\textheight, angle=0]{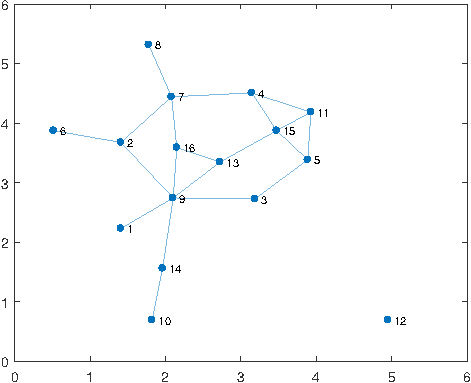}} 
\caption{A graphical representation of the trade network of Renaissance Florentine families, which is considered in Example 2.}
\label{fig.adjacency}
\end{figure}

To test the accuracy of the mean-field approximation proposed
in this work, we make a comparison of the estimates 
using the two methods at hand (Algorithms \ref{alg2} and \ref{alg3}) on a real example taken
form the specialized literature.
The network illustrated in Figure \ref{fig.adjacency} represents the business ties between $16$ Florentine families during Renaissance, and is taken from \cite{Bre1986} (see also \url{http://moreno.ss.uci.edu/data.html}).

The network contains $20$ edges, $47$ 2-stars and $3$ triangles, and is described by the following adjacency matrix:

{
\begin{equation}\nonumber X= \left(  \begin{array}{cccc cccc cccc cccc}
 0 & 0 & 0 & 0&  0&  0&  0&  0&  1&  0 & 0&  0&  0 & 0&  0&  0 \\
 0 & 0&  0&  0&  0&  1 & 1 & 0 & 1 & 0 & 0 & 0 & 0 & 0 & 0 & 0 \\
 0&  0 & 0 & 0 & 1 & 0 & 0 & 0 & 1 & 0 & 0 & 0 & 0 & 0 & 0 & 0 \\
 0&  0 & 0 & 0 & 0 & 0 & 1 & 0 & 0 & 0 & 1 & 0 & 0 & 0 & 1 & 0 \\
 0&  0 & 1 & 0 & 0 & 0 & 0 & 0 & 0 & 0 & 1 & 0 & 0 & 0 & 1 & 0 \\
 0&  1 & 0 & 0 & 0 & 0 & 0 & 0 & 0 & 0 & 0 & 0 & 0 & 0 & 0 & 0 \\
 0 & 1 & 0&  1&  0&  0&  0&  1 & 0 & 0 & 0 & 0 & 0 & 0 & 0 & 1 \\
 0&  0 & 0 & 0 & 0 & 0 & 1 & 0 & 0 & 0 & 0 & 0 & 0 & 0 & 0 & 0 \\
 1& 1 & 1 & 0 & 0 & 0 & 0 & 0 & 0 & 0 & 0 & 0 & 1 & 1 & 0 & 1 \\
 0 & 0 & 0 & 0 & 0 & 0 & 0 & 0 & 0 & 0 & 0 & 0 & 0 & 1 & 0 & 0 \\
 0 & 0 & 0 & 1 & 1 & 0 & 0 & 0 & 0 & 0 & 0 & 0 & 0 & 0 & 1 & 0 \\
 0&  0 & 0 & 0 & 0 & 0 & 0 & 0 & 0 & 0 & 0 & 0 & 0 & 0 & 0 & 0 \\
 0&  0 & 0 & 0 & 0 & 0 & 0 & 0 & 1 & 0 & 0 & 0 & 0 & 0 & 1 & 1 \\
 0&  0 & 0 & 0 & 0 & 0 & 0 & 0 & 1 & 1 & 0 & 0 & 0 & 0 & 0 & 0 \\
 0 & 0 & 0 & 1 & 1 & 0 & 0 & 0 & 0 & 0 & 1 & 0 & 1 & 0 & 0 & 0 \\
 0 & 0 & 0 & 0 & 0 & 0 & 1 & 0 & 1 & 0 & 0 & 0 & 1 & 0 & 0 & 0
 \end{array} \right)\,.
 \end{equation}
}

Also in this case, the mean-field approximation of the moments in the gradient ascent method for maximum log-likelihood estimation is compared with the maximum log-pseudo-likelihood estimate. The proposed method (Algorithm \ref{alg3}) is run for $N_{\rm it}=100000$ iterations and
with $\gamma_{\rm MF}=1 \rm E-4$, and the final estimated parameters are:
$$
\theta_{1, \rm MF}^*=-1.5553 , \quad \theta_{2, \rm MF}^*=-0.0293 , \quad \theta_{3, \rm MF}^*=0.2106.
$$
The time needed for convergence is $165.760$ seconds.
The gradient ascent method for the maximum pseudo-log-likelihood (Algorithm \ref{alg2}) is run for $N_{ \rm it}=10000$ iterations
with $\gamma_{PL}=1 \rm E-3$, and 
the final estimated parameters are:
$$
\theta_{1, \rm PL}^*=-1.6231 , \quad \theta_{2, \rm PL}^*=-0.0188 , \quad \theta_{3, \rm PL}^*=0.2459.
$$
The time needed for convergence is $188.379$ seconds.
The results show that the estimated parameters are in good agreement. However,
the proposed mean-field approximation 
has the advantage of being faster from the computational viewpoint. Likewise in the previous subsection, the final estimated parameters are in the high-temperature phase.

\section{Conclusions}\label{sec:conclusions}

Computational advantages of maximum log-likelihood estimation of a 3-parameters $p$-star model via a modified gradient ascent method based on the mean-field approximation of its 3 moments have been shown, comparing it with maximum log-pseudo-likelihood via gradient ascent. These advantages are evident because, in the first method, empirical quantities (i.e., the empirical moments) are computed only in its initialization. On the contrary, gradient ascent applied to the maximization of the log-pseudo-likelihood function requires the computation of the empirical quantities $\Delta \bX^k_{ij}$ at each gradient step, which is a computationally expensive for a large number of vertices $n$.

The proposed algorithm is applicable whenever the parameters are in the high-temperature phase. This can be easily checked graphically at each iteration (see Figure \ref{fig.AR1}). A significant difference between this algorithm and other approaches for parameters estimation in exponential random graph models is in its use, inside each iteration, of the mean-field approximation of the moments. This is cheap from a computational point of view, since at each iteration one has to solve a nonlinear system with only two equations and two unknowns, even when the number $n$ of vertices in the graph is large. Other approaches, instead, apply at each iteration a more expensive Markov Chain Monte Carlo (MCMC) approximation (such as one based on the Metropolis-Hastings sampler)\footnote{It is known that MCMC estimation procedures do not converge for {\em near degenerate} exponential random graph models \cite{Robinsetal2006}, a condition similar to being in the {\em low-temperature phase} as defined in this paper. The paper \cite{Robinsetal2006} provides also a description of various software packages available for Monte Carlo maximum likelihood estimation, such as {\em SIENA}, {\em pnet}, and {\em statnet}.}. As a possible future extension, the proposed algorithm could be applied to more sophisticated and realistic exponential random graph models, such as the alternating $k$-stars and alternating $k$-triangles models \cite{Snijdersetal2006}, provided the mean-field approximation is still valid in such cases.

\section*{Compliance with Ethical Standards}

\section*{Funding}

The authors acknowledge support from the Italian National Interest Project ``Crisis Lab'' (MIUR, PNR 2011-2013).

\section*{Conflict of Interest}

The authors declare that they have no conflict of interest.

\section*{Appendix}
Approximate values for $p$ and $q$ in the mean-field approximation are found by applying an iterative solver based on Newton's method to the system of Equations \eqref{eq:p}, \eqref{eq:q}:
\begin{align*}
 G_p(p,q):=p- F_p(p,q)&=0, \\
 G_q(p,q):=q-F_q(p,q)&=0.
 \end{align*}
Starting from an initial guess $(p^0,q^0)$, each iteration of Newton's method is given by:
\begin{equation}\label{eq:Newton}
 \left( \begin{array}{c} p^{k+1} \\ q^{k+1} \end{array} \right)=
 \left( \begin{array}{c} p^{k} \\ q^{k} \end{array} \right)+
 \left( \begin{array}{cc} \dfrac{\partial G_p }{\partial p} & \dfrac{\partial G_p }{\partial q} \\
                          \dfrac{\partial G_q }{\partial p} & \dfrac{\partial G_q }{\partial q}   \end{array} \right)^{-1}   \left( \begin{array}{c} p^{k} \\ q^{k} \end{array} \right)
\end{equation}
until $\Vert ( p^k, q^k ) \Vert < {tol}$, for a given tolerance $tol >0$.


\end{document}